\documentstyle{amsppt}

\magnification=1200 \NoBlackBoxes \hsize=11.5cm \vsize=18.0cm

\def\inv{^{-1}}

\def\Ndr{N_d^{\text{red}}}
\def\Ndx{N_d^\times}
\def\Ndt{N_d^{\text{tac}}}

\def\F{\Cal F}

\def\P{\Bbb P}

\def\ss{\vskip.15in}

\def\[{\big[}
\def\]{\big]}
\def\V{\bar{V}}
\def\L{\Cal L}
\def\O{\Cal O}

\def\1/2{\frac{1}{2}}

\def\im{\text{im}}

\def\2{{[2]}}

\topmatter
\title Enumerative geometry of divisorial
families of rational
curves
\endtitle
\author
Ziv Ran
\endauthor

\date May 5, 2002/ Rev. June 10, 2002\enddate

\address University of California, Riverside\endaddress
\email ziv\@math.ucr.edu\endemail
\rightheadtext {Rational curves}
\leftheadtext{Z. Ran}
\abstract We compute the number of irreducible
rational curves of given degree with
1 tacnode or 1 triple point in $\P^2$
or 1 node in $\P^3$
 meeting an appropriate generic
collection of points and lines.
As a byproduct, we also compute the number of
rational plane curves of degree $d$ passing through
$3d-2$ given points and tangent to a given line.
 The method is
'classical', free of Quantum Cohomology.
 \endabstract
 \thanks \raggedright {Updates and corrections
available at math.ucr.edu/$\ \tilde{}\ $
ziv/papers/1nodal.pdf.}\linebreak
Research Partially supported by
NSA Grant MDA904-02-1-0094; reproduction and
distribution by US governement permitted.
\endthanks

\endtopmatter\document

In the past 15 years or so a number of classical
problems in the enumerative geometry of curves in
$\P^n$ were solved, first for $n=2,$ any genus [R1], then
for any $n$, genus 0. The latter development was
initiated by Kontsevich-Manin who developed and used
the rather substantial machinery of Quantum Cohomology
(cf. e.g. [FP]).
Subsequently, in a series of papers [R2-R5] the author
developed an elementary alternative approach,
free of Quantum
cohomology, and used it to solve a number of classical
enumerative problems for rational, and sometimes elliptic,
curves in $\P^n, n\geq 2.$ The present paper continues
this series. The object here is to enumerate the
irreducible rational curves of given degree $d$
in $\P^2$ with one tacnode or one triple point
passing through
$3d-2$ general points (see Theorems 2,5 below),
as well as the irreducible
rational curves
in $\P^3$ with one ordinary node which
contain $a$
general points
and are incident
to $4d-2a-1$ general lines
(see Thm 1 below).
Note that the family of
1-tacnodal (resp. 1-triple point, resp. 1-nodal)
curves in $\P^2$
(resp. $\P^2, \P^3$) is of codimension
1 in the family of all rational curves so that we are
effectively computing the 'degree', in a sense,
of certain natural divisors in the family of all
rational curves.
Indeed by a result of Diaz and Harris [DH]
in the case of $\P^2$, the general
member of any such divisor, if not nodal, is either
1-cuspidal (which case was enumerated in [R2])
or 1-tacnodal or has 1 triple point.
It seems very likely, but doesn't seem to be
in the literature, that the natural analogue of this
result also holds for any $n\geq 3$: i.e. that
the general member of any divisor in the family of
rational curves is either smooth, or 1-nodal reducible
(any $n$), or 1-nodal irreducible (n=3).
As a byproduct of the proof of Theorem 2,
we also obtain a formula for one of the
'characteristic numbers' for rational plane curves,
viz. the number of such curves of degree $d$ passing
through $3d-2$ general points and tangent to a given line
(cf. Cor. 4 below).
\par
Enumeration of rational
plane curves with a tacnode or triple point
appears to be well within the range of
interest, at least, of classical geometers.
For quartics, results of this kind, at least
the enumeration of quartics with a triple point,
were obtained by Zeuthen in 1882. Recently,
enumerative results of this kind for curves with
'few' singularities on general surfaces were obtained
by Kleiman and Piene \cite{KP}.
Enumeration of 1-nodal rational curves in
$\P^n$ for any $n$ was recently
announced by Zinger \cite{Z}, using Quantum Cohomology.
As we shall see below, the 1-nodal case in
$\P^3$
case is analogous to, but
easier than the 1-tacnodal and 1 triple-point
cases in $\P^2$. Our proof is based on the intersection
calculus
on the nonsingular model of the surface swept out by
the appropriate 1-parameter family of rational curves,
developed in earlier papers [R2-R4],
together with elementary
residual-intersection (in particular,
double-point) theory as in [F]

We begin by reviewing some qualitative
results about families
of rational curves in $\P^n$, especially for
$n=2$ or $3$.
See [R2][R3] [R4] and references
therein for details and proofs.
In what follows we denote by $\bar{V}_d$ the
closure in the Chow variety of the locus of
irreducible
nonsingular rational
curves of degree $d$ in $\P^n, n=3$,
with the scheme structure
as closure, i.e. the reduced structure
(recall that the Chow form
of a reduced 1-cycle $Z$ is just
the hypersurface in $G(1,\P^3)$ consisting
of all linear spaces meeting $Z$); if $n=2$
$\bar{V}_d$ is just the (closure of)
the Severi variety. Thus
$\V_d$ is irreducible reduced of dimension
$$\dim (\V_d)=4d, n=3$$
$$\dim (\V_d)=3d-1, n=2.$$
Let $$A_1, \ldots, A_k\subset\P^n$$ be a
generic collection of linear subspaces of respective
codimensions $a_1,
\ldots, a_k,$ $ 2\leq a_i \leq n$
(so if $n=2$ these are just points).
We denote by
$$
B= B_d = B_d (a_{\cdot}) = B_d (A_{\cdot})
$$
the normalization of the locus (with reduced structure)
$$
\{(C, P_1, \ldots, P_k ) \ :
\ C \in \bar{V}_d , P_i  \in C \cap A_i, i =1
, \ldots, k \}
$$
which
is also the normalization of its
projection to $\bar{V}_d$, i.e. the
locus of degree-$d$ rational curves
(and their specializations) meeting
$A_1, \ldots, A_k$.
 We have
$$
\dim B = (n+1)d+n-3  - \sum (a_i - 1) .\tag 1
$$
When $\dim B=0$ we set
$$N_d(a.)=\deg (B). \tag 2$$
When $n=2$, so all the $a_i=2$ they will be dropped.
The integer
$k$ is called the {\it length} of the
condition-vector $(a.)$.
The numbers $N_d$ and $N_d(a.)$, first
computed in general by Kontsevich and Manin
(see for instance
[FP] and references therein),
 are computed in [R2],[R3] by an elementary method
based on recursion on $d$ and $k$.

Now suppose $\dim B=1$ and let
$$
\pi : X \to B \tag 3
$$
be the normalization of the tautological family of rational
curves, and $$f: X \to \P^n$$ the natural map.
The following summarizes results from [R2][R3][R4] :
\proclaim{Theorem 0}(i) X is smooth .\par
(ii) Each fibre $C$ of $\pi$ is either\par
(a) a $\P^1$ on which $f$ is
either an immersion with at most one exception
which maps to a cusp
($n=2$) or an embedding ($n>2$); or\par
(b)  a pair of $\P^1$'s meeting transversely once,
on which $f$ is an
immersion with nodal image ($n=2$)
or an embedding ($n>2$); or\par
(c)  if $n=3$, a $\P^1$ on which $f$ is a
degree-1 immersion such that $f(\P^1)$ has
a unique singular point which is an ordinary node.\par
(iii) If $n>2$ then $\bar{V}_{d,n}$
is smooth along the image $\bar{B}$ of
$B$, and $\bar{B}$ is smooth except,
in case some $a_i=2,$ for
ordinary nodes corresponding to curves meeting some
$A_i$ of codimension 2 twice.
If $n=2$ then $\bar{V}_{d,n}$ is smooth
in codimension 1 except for a cusp along
the cuspidal locus and
normal crossings along the
reducible locus, and $\bar{B}$ has the
singularities induced from
$\bar{V}_{d,n}$ plus ordinary nodes
corresponding to curves with a node at
some $A_i$, and no other
singularities.
\endproclaim

Next, we review some of the enumerative apparatus introduced
in [R3][R4] to study $X/B$.
Set
$$m_i=m_i(a.)=-s_i^2, i=1,...,k.\tag 4 $$
Note that if $a_i=a_j$ then $m_i=m_j;$
in particular for $n=2$ they are all equal
and will be denoted by $m_d.$
It is shown in [R2] [R3][R4] that these numbers can all be
computed recursively in terms of data of lower degree $d$ and
lower length $k$.
For instance for $n=2$ we have
$$
2m_d = \sum_{d_1 + d_2 = d}
N_{d_1} N_{d_2} d_1 d_2 {\binom{3d-4}{3d_1-2}}.\tag 5
$$
For $n>2,$ note that
$$
s_i . s_j = N_d(...,a_i+a_j,...,\hat a_j,...), i \neq j.
$$

(so for $n=2$ this is always 0).
Also, letting $R_i $ denote the sum of all
fibre components not meeting $s_i$ , we have
$$
s_1 \cdot R_2 = \sum N_{d_1} (A_{\cdot}^1 , A_1,\P^{s_1} )
N_{d_2} (A^2_{\cdot},
A_2, \P^{s_2} ).
\tag 6
$$
the summations being over all $d_1 + d_2 = d, s_1 + s_2 = 3 $
and all decompositions
$A_{\cdot} =
 (A_1, A_2) \coprod (A_{\cdot}^1) \coprod (A_{\cdot}^2)
$
(as unordered sequences or partitions);
similarly for the other
$s_i.R_j$.
So all these numbers may be considered known.
Then we have
$$m_i=\frac{1}{2}(s_i.R_j+s_i.R_p-s_j.R_p)-s_i.s_j
-s_i.s_p+s_j.s_p $$
for any distinct $i,j,p,$ and the RHS here is an expression of
lower degree and/or length, hence may be considered known.
\par
Next, set $$L=f^*(\O(1)),$$ and note that
$$L^2=N_d(2,a.),\  L.s_i=N_d(a_1,...a_{i+1}...)
,\ i=1,...,k$$
(in particular, $L.s_i=0$ if $a_i=n.$)
We computed in [R3] that, for any $i$,
$$
L \sim d s_i - \sum\limits_{F \in \F_{i}} \deg (F) F +
( N_d (a_1 ,\ldots, a_i+1, \ldots) + dm_i) F_0
$$
where $F_0$ is the class of a complete fibre
and $\F_i$ is the set of fibre components not meeting $s_i$.
Consequently we have
$$
N_d(2,a_1, \ldots) =
2d N_d (a_1 + 1, a_2, \ldots) + d^2 m_1(a.) -
\sum\limits_{F\in \F_{1}(a.)} (\deg F)^2
\tag 7
$$
and clearly the RHS is a lower degree/length expression, so
all the $N_d(2,\ldots)$ are known. We also have
for $n>2$ that
$$
N_d(a_1,a_2+1,...) - N_d(a_1+1,a_2,...)=$$$$
dN_d(a_1+a_2,...)-\sum\limits_{F
\in (\F_{1}-\F_{2})(a.)} (\deg F) + N_d(a_1+1,a_2,...)+dm_1(a.)
\tag 8
$$
and again the RHS here is 'known', hence so is the LHS,
which allows us to 'shift weight' between the $a_i$'s
till one of them becomes 2,
so we may apply (7), and thus
compute all of the $N_d(a.)$'s.\par
Next, it is easy to see as in [R3] that
$$L.R_i=
\sum_{d_1+d_2=d}\binom{3d-1}{3d_1-1}d_1d_2^2N_{d_1}N_{d_2},
\  n=2\tag 9$$
(in this case this is independent of $i$
and we will just write it as $L.R$);

$$
L.R_i = \sum d_2N_{d_1} (a_{\cdot}^1 ,{s_1} )
N_{d_2} (a^2_{\cdot},
 {s_2} ),\ n>2
\tag 10 $$
the summation for $n>2$ being over all
$d_1 + d_2 = d, s_1 + s_2 = 3 $
and all decompositions
$$A_{\cdot} =  (A_{\cdot}^1) \coprod (A_{\cdot}^2)
$$ (as unordered sequences or partitions) such that
$A_i\in(A_{\cdot}^1)$.\par
Finally, the relative canonical class
$K_{X/B}=K_X-\pi^*(K_B)$ was computed in [R3] as
$$K_{X/B}=-2s_i-m_iF+R_i\tag 11$$
for any $i$. Note that $-R_i^2$ equals
the number of reducible fibres in the family $X/B$,
a number we denote by $N_d^{\text{red}}(a_.)$, and
which is easily computable by recursion.
From this we compute easily that
$$L.K_{X/B}= -2N_d(...a_i+1...)-dm_i+L.R_i,\tag 12$$
$$K_{X/B}^2=-\Ndr(a_.).\tag 13$$
Now in case $n=3$,
for any condition-vector $(a_.)$ of weight
$$\sum (a_i-1)=4d-1,$$
we denote by $\Ndx(a_.)$ the number of 1-nodal
irreducible rational curves in $P^3$ meeting
a generic collection of linear spaces of
respective codimensions $(a_1,...,a_k)$.
\proclaim{Theorem 1} We have for any $i=1,...,k,$
\ss\noindent
$(14)\ \ \Ndx(a_.)=$\ss\noindent
$$(d-3)N_d(2,a_.)-N_d(...a_i+1...)
-(4d+2)m_i+2L.R_i+\Ndr(a_.).$$\endproclaim
\demo{proof}
We set things up so as to apply the {\it double point
formula} (in a relative form [F]).
Let $$X_B^2$$ denote the fibre square
of $X/B$, and
$$\Delta\simeq X\subset X_B^2$$ the
diagonal. Then $X_B^2$ is smooth except at points
$(p,p)$, where $p$ is a singular point of
a fibre of $X/B$. Locally at such points,
if $X/B$ is given locally at $p$ by
$$xy=t,$$
then $X_B^2$ is given by
$$x_1y_1=x_2y_2=t,$$
and therefore has an ordinary 3-fold double point.
Moreover the diagonal $\Delta$ is defined by
$$x_1=x_2, y_1=y_2$$
and in particular is a non-Cartier divisor.
Let $$b: Y\to X_B^2\tag 15$$
denote the blow-up of  $\Delta$, with
exceptional divisor $$\Delta'=b^*(\Delta).$$
Then it is easy to see that $b$ is a small
resolution of $X_B^2$ with exceptional locus
$$E=\sum E_p$$
consisting of a $\P^1$ for each relatively
singular point $p$. Moreover the fibre of $Y$
over $B$ corresponding to a reducible fibre
$C_1\cup_pC_2$ of $X/B$ is of a 'honeycomb'
shape
$$\matrix B_{(p,p)}C_1^2&\cup&C_1\times C_2\\
C_2\times C_1&\cup&B_{(p,p)}C_2^2\endmatrix
$$
where
$$B_{(p,p)}C_1^2\cap B_{(p,p)}C_2^2=E_p,
C_1\times C_2\cap C_2\times C_1=\emptyset.$$
Also $\Delta'\to\Delta$ is just the blowing up
of all the points $p$. Furthermore, the
identity of rational functions, locally at $p$,
$$\frac{y_2-y_1}{x_2-x_1}=\frac{y_1}{x_2}=
\frac{y_2}{x_1}$$
shows that there $b:Y\to X_B^2$ coincides
with the blowup of the ideal $(y_1,x_2)$,
and with that of the ideal $(y_2,x_1)$.\par
We will need to know the normal bundle
$$\nu=N_{\Delta'/Y}.$$
To this end, note that, clearly
$$\omega_{\Delta'}=b^*\omega_X(E).$$
On the other hand, the map
(15) is small, hence crepant, so
$$\omega_Y=b^*\omega_{X^2_B}=b_1^*\omega_X\otimes
b_2^*\omega_X\otimes c^*\omega_B\inv\tag 16$$
where we use the evident maps
$$b_i:Y\to X, i=1,2, c:Y\to B.$$
Therefore by the adjunction formula we conclude
$$\nu=b^*\omega_{X/B}\inv(E).\tag 17$$
Now we are ready to apply the double-point formula
to the map
$$f^2=(f_1,f_2):X_B^2\to \P^3\times\P^3.$$
This shows that the cycle of ordered pairs $(x_1,x_2)$
in $X/B^2$ such that $f(x_1)=f(x_2)$
is residual to $\Delta'$ in $f^{2*}(\Delta_{\P^3})$.
Therefore,
as in ([CR], Thm 3) we find that
$$\Ndx(a_.)=\frac{1}{2}c_3(f_1^*(\O(1))\otimes
f_2^*(Q)\otimes
\O(-\Delta')),$$
where $Q$ denotes the universal quotient bundle on
$\P^3.$ Then a straightforward calculation,
based on the intersection calculus reviewed above,
 yields
the formula (14).\qed\enddemo
We turn next to the enumeration of rational plane
curves with a tacnode. Denote by $\Ndt$
(resp. $\kappa_d$) the number of
rational curves with a tacnode (resp. cusp) passing through
$3d-2$ generic points in $\P^2$, and recall that $\kappa_d$
was computed in [R3].
\proclaim{Theorem 2} We have
$$\Ndt=\hskip7cm \tag 18$$
$$\hskip2cm (d-4)N_d-(d-7)(d^2+7d-4)m+\frac{d-7}{2}L.R
+2\Ndr+\frac{1}{2}\kappa_d $$
\endproclaim
\demo{proof} We use an analogous setup and
notation as above, this time
for a pencil of rational curves in $\P^2$ through $3d-2$
points. To begin with, we describe an algebraic
setup for the (relative) Gauss mapping. Consider the
relative principal parts sheaf for $L$ on $X/B$, which
fits in an exact sequence
$$0\to \Omega_{X/B}(L)\to P_{X/B}(L)\to L\to 0.\tag 19$$
Setting $V=H^0(\O_{\P^2}(1))$, we have a natural
map $V\to P_{X/B}(L)$ which combines with the Euler
sequence to yield an exact diagram
$$\matrix 0&\to& f^*\Omega_{\P^2}(1)&\to& V\otimes\O_X&
\to&f^*\O_{\P^2}(1)&\to&0\\
&&\downarrow &&\downarrow&&\parallel&&\\
0&\to&\Omega_{X/B}(L)&\to&P_{X/B}(L)&\to&L&\to&0.\endmatrix$$
Let $P'$ denote the pushout of $P_{X/B}(L)$ by the
natural map
$$\Omega_{X/B}(L)\to \omega_{X/B}(L).$$
Now let $b:X'/B\to X/B$ denote the blowup of all singular
points of fibres and all 'cuspidal' points, i.e. all points
$(x,b)\in X$ such that $f(x)$ is a cusp on $f(X_b)$ (see
[R3] for a discussion and enumeration of these).
Let $E=\sum E_p$ be the exceptional divisor over the
singular points of fibres and $U=\sum U_q$ be the exceptional
divisor over the cuspidal points.
Set $f'=b\circ f$.
Then it is easy to see by a local computation that
\proclaim {Lemma 3}
The image of the natural map
$$f^{'*}\Omega_{\P^2}\to b^*\omega_{X/B}$$
coincides with $b^*\omega_{X/B}(-E-U)$
\endproclaim\demo{proof}
It suffices to prove locally at each fibre node
or cuspidal point $p$
that
$$\im(f^*\Omega_{\P^2}\to \omega_{X/B})=\omega_{X/B}.I_p.$$

In local coordinates at $p$, the family $X/B$ is given by
$$xy=t$$
with $f=(x,y)$
and $\omega_{X/B}$ is generated by $dx\wedge dy.dt\inv.$
Since
$$dx\wedge dt = xdx\wedge dy, dy\wedge dt =-ydx\wedge dy,$$
we have that $f^*\Omega_{\P^2}$ is generated by
$xdx\wedge dy.dt\inv, ydx\wedge dy.dt\inv,$ as claimed.\par
In the cuspidal case, our family has local
coordinates $u,t$ with
$$f=(u^2,(u^2-t)u)=(x,y).$$
Here $\omega_{X/B}$ is generated by $du$ and its subsheaf
generated by $$dx=2udu, dy=(3u^2-t)du$$ clearly coincides
with the subsheaf generated by
$udu$ and $tdu,$ as claimed\qed\enddemo

It follows from Lemma 3 that if we let
$P$ denote the image of the natural map
$$V\otimes\O_{X'}\to b^*P',$$ then we get a diagram
with exact rows and columns
$$\matrix 0&\to& f^{'*}\Omega_{\P^2}(1)&\to& V\otimes\O_X&
\to&f^*\O_{\P^2}(1)&\to&0\\
&&\downarrow&&\downarrow&&\downarrow&\\
0&\to&b^*\omega_{X/B}(L-E-U)&\to&P&\to&b^*L&\to&0\\
&&\downarrow&&\downarrow&&\downarrow&\\
&&0&&0&&0&\endmatrix.\tag 19$$
Now let
$$\phi:I\to \P^2$$
be the incidence (or flag) variety, with universal
flag
$$V\otimes\O_I=F^0\supset F^1=\phi^*\Omega_{\P^2}(1)
\supset F^2\supset (0),$$
$$F^0/F^1=\phi^*\O(1).$$
Now the diagram (19) gives rise to a lifting of $f'$
to a morphism
$$g:X'\to I$$
with
$$g^*(F^0/F^1)=L':=b^*(L),
g^*(F^1/F^2)=b^*(\omega_{X/B})(- E- U).
\tag 20$$
Clearly, the value of $g$ at a point
$$(x,b)\in X'\setminus E\setminus U$$
mapping to $$y\in\P^2$$ is the pair
$$(y,T_y(f(X_b))$$
consisting of $y$ and the tangent line to
$f(X_b)$ at $f(y)$; in particular, the tacnodes
in the family $X/B$ correspond (1:2) to the double points of
$g$, which we propose to count as in the proof of Theorem 1.
To this end, the construction of a 'good'
desingularization of $X'\times_BX'$ must be modified
to take into account the fact that each $E_p$ appears
in its fibre over $B$ with multiplicity 2.\par
Consider then the fibre square
$$\bar{Y}=X'\times_B X'.$$
Its singular points are as follows. First, for each
cuspidal point $q$ on $X$, we get a 3-fold ODP $(y,y)$
where $y$ is the intersection of $U_q$ with
the proper transform of the fibre of $X/B$ containing $q$.
We desingularize this by blowing up $U_q\times U_q$
(which is the same locally as blowing up the diagonal).
Here the construction is exactly as in the proof of Thm 1.
\par
Additionally, $\bar{Y}$ is singular along each $E_p\times
E_p$, where it has a local equation of the form
$$x_1^2u_1=x_2^2u_2$$
where $x=0$ is the equation of $E_p$, $u$ is a coordinate
on $E$ and $x_i,u_i$ are the pullbacks of $x,u.$
This can be desingularized by blowing up the locus
$$x_1=x_2=0,$$
i.e. $E_p\times E_p.$ The exceptional divisor maps
to $E_p\times E_p$ generically with degree 2,
ramified over the
4 rulings
$$E_p\times\{0,\infty\}\cup \{0,\infty\}\times E_p,$$
where $0,\infty$ are the 2 intersections of $E_p$ with
the other 2 components of its fibre, and having
$\P^1$ fibres over $\{0,\infty\}^2$.\par
Let
$$b: Y\to \bar{Y}$$
be the global desingularization thus obtained, and
$\Delta"\subset Y$ the proper (=total) transform
of the diagonal. It is easy to see that the map
$$b:\Delta"\to X'$$ identifies $\Delta"$ with the
blowup $X"$ of all singular points of set-theoretic fibres
of $X/B$
(i.e. all the points of the form $(y,y)$ or $(0,0)$
or $(\infty,\infty)$ in the above notation.
\par
The normal bundle
$$\nu=N_{\Delta"/Y}$$
can be computed much like before: first,
$$\omega_{\Delta"}=b^*\omega_{X'}(E")$$
where $E"=\sum E"_r$ is the exceptional divisor of
$b|_{\Delta"}.$ Next, as before,
$$\omega_{\bar{Y}}=b_1^*\omega_X\otimes b_2^*\omega_X\otimes
\omega_B\inv.$$
As $\bar{Y}$ has an ordinary double surface generically
along
$\bigcup_pE_p$, it follows that
$$\omega_Y=b^*(\omega_{\bar{Y}})(-\bigcup_pE_p).$$
Putting these together we see that, using
divisor notation and setting $K=K_X-K_B,$
$$\nu=-(b^*(\omega_{X'}-\omega_B-E)-E")=-(K+U-E").\tag 21$$
Consequently, we have
$$\nu^2=K^2-2\Ndr-2\kappa_d\tag 22$$
where $\kappa_d$ is the number of cuspidal
rational plane curves of degree $d$ through $3d-2$
generic points, computed in [R3]. Of course,
we also have
$$L.\nu = -L.K = 2dm-L.R\tag 23$$
\par
Now, to use double-point theory
 consider the cartesian product
$I^2=I\times I$ with projections
$$p_i:I^2\to I, i=1,2, \phi^2:I^2\to(\P^2)^2.$$
Let $F_.$ be the dual filtration to $F^.$,
defined by
$$F_i=(F^0/F^{i})^*\subset (F^0)^*.$$
In $I^2,$ $\phi^{2*}\Delta_{\P^2}$ is defined as the zero-
locus of a map
$$p_1^*F_1\to p_2^*(F_3/F_1),$$
and inside this, $\Delta_I$ is defined as the zero-
locus of a map
$$p_1^*(F_2/F_1)\to   p_2^*(F_3/F_2).$$
Considering as before the map
$$f^2=(f_1,f_2):Y\to I^2,$$
double-point (or more precisely,
residual-intersection) theory shows that\ss\noindent
$2\Ndt =$\ss\noindent
$$c_2(f_1^*(F_1^*)\otimes f_2^*((V^*\otimes\O)/F_1)
\otimes\O(-\Delta"))c_1(f_1^*((F_2/F_1)^*)\otimes
f_2^*(F_3/F_2)\otimes\O(-\Delta")).$$\ss\noindent
Using (20), the $c_2$ factor can be identified,
writing $L_i=p_i^*L, K_i=p_i^*(\omega_{X/B})$,  as
$$L_2^2+L_1L_2-L_2\Delta'+(L_1-\Delta')^2$$
while the $c_1$ factor is
$$(K_1-E-U)+L_2-\Delta".$$
Then a routine computation, using as above the
necessary intersection theory on $Y$ and $X$
yields the formula (18).\qed

\enddemo
As a bonus for the construction of the lift $g$
above- and without using any double-point theory,
we get a formula for one of the 'characteristic
numbers', namely the number $N_d^t$ of
rational plane curves through $3d-2$ points tangent
to a given line. Thus let
$$p_2:I\to\P^{2v}$$
be the natural map of the incidence variety to the dual
projective plane, and set
$$f^v=p_2\circ g:X'\to\P^{2v}, L^v=f^{v*}\O(1).$$
As above, we can write
$$L^v\sim (2d-2)s_1-\sum
\limits_{F\in\Cal F_1} (\deg(F)-1)F'-\sum E_p
-\sum U_q+xF_0$$
where the first sum is over all fibre components
of $X/B$ not meeting $s_1$, $F'$ denotes the proper
transform of $F$ on $X'$, the second and third sums
cover over all
exceptional divisors of $X'/X$,
$F_0$ denotes a fibre of $X'/B$, and we have simply used
the fact that the dual to a rational curve of degree $d$
has degree $2d-2$. Now intersect this expression with $s_1$
and note that
$$L^v.s_1=N_{d\to},$$
i.e. the number of degree-$d$ rational curves
through $A_1,...,A_{3d-2}$ with a given tangent
direction at $A_1.$ This yields
$$x=N_{d\to}+(2d-2)m_d.$$
Since $N_d^t=\L^{v2},$ this number can now be easily
computed from our intersection calculus, yielding
the following
\proclaim{Corollary 3} We have
$$N_d^t=4(d-1)^2m_d+d(d-1)N_{d\to}+\Ndr
-\kappa_d-2\sum\limits
_{F\cap s_1=\emptyset}(\deg(F)-1)^2.$$\endproclaim

Finally, we take up the problem of enumerating the
triple points, i.e. of computing
the number $N_d^{\text{tri}}$ of irreducible rational
curves with a triple point passing through $3d-2$
generic points. Again our tool will be a suitable
residual-intersection computation. Using notations as
above, take 3 copies $$Y_{12}, Y_{13}, Y_{23}$$
of the small resolution of $X_B^2$
constructed above. We have evident 'first projection'
maps
$$Y_{12}\to X, Y_{13}\to X$$
which allows us to construct
$$Y_{12}\times_XY_{13}$$
together with a map
$$p_{23}:Y_{12}\times_XY_{13}\to X_B^2.$$
Let $Z$ be the unique component of
$Y_{12}\times_XY_{13}$ which dominates $ X_B^2.$
Define $$Y^*$$
as the blowup of $Z$ in the proper transform of the
diagonal in $X^2_B.$
Then $Y^*$ comes equipped with projections
$$p_{ij}:Y^*\to Y_{ij},$$
a map
$$f^3=(f_1,f_2,f_3):Y^*\to (P^2)^3$$
as well as with a birational map
$$b^*:Y^*\to X^3_B.$$
An elementary analysis shows that $Y^*$
has exactly 2 singular points, each
isomorphic to a cone over a Segre variety
$\P^1\times\P^2\subset\P^5$, hence also to
the germ at the origin of the determinantal
variety $M^1_{2\times 3}$ of $(2\times 3)$
matrices of rank at most 1, and
consequently each of these points admits a small
'determinantal' resolution with $\P^1$ as exceptional
locus.
It can be shown that this resolution admits a natural
morphism to the relative Hilbert scheme
$\Cal{H}$ilb$_3(X/B).$
Though we could replace $Y^*$ by this resolution,
this will not be necessary since by definition
the big diagonals are already Cartier on
$Y^*$.\par
Now consider on $Y^*$ the locus
$$D_{12}=p_{12}^*D$$
where $D\subset Y$ is the double locus of $f^2$,
i.e. the closure of the locus of distinct point
$(y_1,y_2)$ in the same fibre which map to a node,
i.e. such that $f(y_1)=f(y_2)$. Then clearly
$$[D_{12}]=p_{12}^*([D])=c_2(Q_2(L_1-\Delta_{12}))$$
where
$$Q_i=f_i^*Q, L_i=f_i^*L, \Delta_{12}=p_{12}^*(\Delta').$$
Now consider the intersection $D_{12}.D_{13}$.
This consists of the locus we want, i.e. that
of ordered distinct triples $(y_1,y_2,y_3)$ on the same
fibre such that $f(y_1)=f(y_2)=f(y_3)$, plus the locus
$D_{12}.\Delta_{23}$. Then residual-intersection theory
tells us that
$$6N_d^{\text{tri}}=c_2(Q_2(L_1-\Delta_{12}))
.c_2(Q_3(L_1-\Delta_{13}-\Delta_{23})).$$
This   expression can be computed by routine
calculations, using the interesection calculus developed
above; the only possibly nonobvious terms are
$$\Delta_{12}^2\Delta_{13}^2=\Delta_{12}^2\Delta_{23}^2
=K^2,$$
$$\Delta_{12}^2\Delta_{13}\Delta_{23}=(K-E)^2.$$
The former follows from the fact that
$$p_{12*}(\Delta_{13}^2))=-p_{12*}(p_{13}^*(K_1-E))=-K_1,$$
and that $\Delta_{12}^2=-(K_1-E)$ (cf. (17)).
The latter follows from the fact that
$$p_{12*}(\Delta_{13}\Delta_{23})=\Delta_{12}.$$
 This
yields
\proclaim{Theorem 5} We have\ss\noindent
$6N_d^{\text{tri}}=$\ss\noindent
$$(3d^2-18d+30)N_d+(3d-18)(dm_d-L.R)-6\Ndr$$
\endproclaim
\remark{Check} For $d=4$, we have $m_d=428, N_d=620,
L.R=3276, \Ndr=2124$, therefore
$$N_4^{\text{tri}}=60,$$
a number computed classically by Zeuthen
and recently, with modern methods, by Kleiman and Piene
\cite{KP} .
I am grateful to Steve Kleiman for this information,
which led to a correction of an error in an
earlier statement of the formula.\endremark


\Refs\widestnumber \key{CR} \ref\key CR\paper Dimension of
families of space curves \by M.C.Chang and Z. Ran\jour Comp.
Math.\vol 90 \pages 53-57\yr 1994\endref
\ref\key DH\paper
Geometry of Severi varieties
\by S. Diaz and J. Harris\jour Trans. Amer. Math.
Soc.\vol 309\yr 1988\pages 1-34\endref
\ref\key F\by W.
Fulton\book Intersection theory \publ Springer\yr 1984\endref
\ref\key FP\bysame and R. Pandharipande
\paper Notes on stable maps
and quantum cohomology \inbook Algebraic Geometry
\procinfo Santa
Cruz 1995\endref
\ref\key KP\by S. Kleiman, R. Piene
\paper Enumerating singular curves on surfaces
\inbook Algebraic geometry,
Hirzebruch 70, Contemporary Math.\vol 241
\yr1999\pages 209-238
\finalinfo corrections and
revision in math.AG/9903192 \endref
 \ref\key R1\by Z.Ran\paper
Enumerative geometry of singular plane curves\jour Invent.
math\vol 97\yr 1987\pages 447-465
\endref \ref\key R2\bysame\paper Bend, break and count \jour
Isr. J. Math\vol 111 \yr 1999 \pages109-124\endref \ref\key
R3\bysame \paper Bend, break and count II
\jour Math. Proc. Camb.
Phil . Soc. \vol 127\yr 1999\pages 7-12
\endref
\ref\key R4\bysame\paper On the variety of rational space
curve\jour Isr. J. Math \vol 122\yr 2001\pages 359-370\endref
\ref\key R5\bysame\paper The degree of the divisor of jumping
rational curves\jour Quart. J. Math.\yr 2001
 \pages 1-18\endref
 \ref\key
Z\by A. Zinger\paper Enumeration of 1-nodal rational curves in
projective spaces\jour math.AG\vol 0204236
\endref

\endRefs
\enddocument